\documentclass[12pt]{amsart}
\newif\ifdraft\draftfalse

\usepackage{amssymb}
\usepackage{amsthm}
\usepackage{eucal}
\usepackage[english]{babel}

\usepackage{parskip}

\usepackage{tikz,pgfplots}

\usepackage{rotating}

\usepackage{tikz-cd}

\usetikzlibrary{cd}

\usetikzlibrary{positioning}
\usetikzlibrary{arrows}
\usepackage{xcolor}

\usepackage{mathrsfs} 

\usepackage{manfnt}
\usepackage{color}

\makeatletter
\def\@begintheorem#1#2[#3]{%
    \def\naam{#1}
  \deferred@thm@head{\the\thm@headfont \thm@indent
    \@ifempty{#1}{\let\thmname\@gobble}{\let\thmname\@iden}%
    \@ifempty{#2}{\let\thmnumber\@gobble}{\let\thmnumber\@iden}%
    \@ifempty{#3}{\let\thmnote\@gobble}{\let\thmnote\@iden}%
    \thm@swap\swappedhead\thmhead{#1}{#2}{#3}%
    \the\thm@headpunct
    \thmheadnl 
    \hskip\thm@headsep
  }%
  \ignorespaces}
\makeatother

\newcommand{\kantlijndraft}[1]{\ifdraft\hspace{-\lastskip}%
\vadjust{\vspace{-1mm}\smash{\llap{{\tt #1}\hspace{8mm}}}\vspace{1mm}}\fi}

\def\voegToe#1#2#3{\immediate\write1{\string\newlabel{#1}{{#2}{#3}}}}

\newcommand{\thlabel}[1]{\voegToe{#1}{\naam\noexpand~\thetheorem}{\thepage}\kantlijndraft{#1}}

\makeatletter
\renewcommand{\label}[1]{\voegToe{#1}{\@currentlabel}{\thepage}\kantlijndraft{#1}}
\makeatother

\newtheorem{theorem}{Theorem}
\newtheorem{lemma}{Lemma}

\newtheorem{ECT}{Extension Construction Technique}
\newtheorem{proposition}{Proposition}
\theoremstyle{definition}
\newtheorem{definition}{Definition}
\newtheorem{example}{Example}

\newtheorem{comment}{Remark}
\theoremstyle{remark}

\newtheorem{claim2}{\sc Claim}



\title{$n$-H-closed Spaces}
\date{}
\author{Fortunata Aurora Basile}
\address[Basile]{University of Messina}
\email{basilef@unime.it}
\author{Maddalena Bonanzinga}
\address[Bonanzinga]{University of Messina}
\email{mbonanzinga@unime.it}
\author{Nathan Carlson}
\address[Carlson]{California Lutheran University}
\email{ncarlson@callutheran.edu}
\author{Jack Porter}
\address[Porter]{University of Kansas}
\email{porter@ku.edu}

\begin{document}
\maketitle

\begin{abstract}
In this paper we extend the theory of H-closed extensions of Hausdorff spaces to a class of non-Hausdorff spaces,  defined in \cite{B}, called $n$-Hausdorff spaces. The notion of H-closed is generalized to an $n$-H-closed space. Known construction for Hausdorff spaces $X$, such as the Kat\v{e}tov H-closed extension $\kappa X$, are generalized to a maximal $n$-H-closed extension denoted by $n$-$\kappa X$.
\end{abstract}

{\bf Keywords: }$n$-Hausdorff spaces; H-closed spaces; Kat\v{e}tov H-closed extension.

{\bf AMS Subject Classification:} 	54D10, 54D20, 	54D35, 	54D80.

\section{Introduction}

\hspace{\parindent} A Hausdorff space that is closed in every Hausdorff space in which it is embedded is called {\bf H-closed} - short for Hausdorff-closed. H-closed spaces were introduced in 1924 by Alexandroff and Urysohn \cite{AU} and characterized as  those Hausdorff spaces in which every open cover has a finite subfamily whose union is dense.  It is straightforward to show that  every open cover of a space $X$ has a finite subfamily whose union is dense iff for every open filter $\mathcal{F}$ on $X$, $a\mathcal{F} \ne \varnothing$ iff for every open ultrafilter $\mathcal{U}$ on $X$, $a\mathcal{U} \ne \varnothing$ (Hausdorff is not needed). Most of the basic properties of H-closed spaces appear in  \cite{PW}.

\hspace{\parindent} The theory of H-closed extensions of a non H-closed Hausdorff space $X$ is well developed (see \cite{PW}). The Kat\v{e}tov H-closed extension $\kappa X$ defined on the set $X\cup \{{\mathcal U}: {\mathcal U}$ is a free open ultrafilter on $X\}$ is the projective maximum among all H-closed extensions of $X$. The Fomin extension $\sigma X$ is $\sigma X=(\kappa X)^\#$, where $Y^\#$ is the strict extension defined in \S 3, is another important $H$-closed extension of $X$.  

\hspace{\parindent} In 2013, Bonanzinga \cite{B} expanded the property  of Hausdorff to n-Hausdorff  where $n \in \omega,\;n\geq2$ by defining a space to be  {\bf n-Hausdorff} if for each distinct points $x_1,...,x_n\in X$, there exist open subsets $U_i$ of $X$ containing $x_i$ for every $i=1,...,n$ such that $x_i \in U_i$ and $\bigcap_{i=1}^{n}U_i = \varnothing$.  Clearly, we have Hausdorff when $n = 2$. 

\hspace{\parindent} In this study, we develop the theory of $n$-H-closed spaces, for every $n\in\omega$, $n\geq 2$ that generalizes the theory of H-closed spaces. This theory is developed within the context of $n$-Hausdorff spaces, as H-closed spaces are studied within the class of Hausdorff spaces. A maximal $n$-H-closed extension $n$-$\kappa X$ is constructed as well as the related Fomin extension $n$-$\sigma X=(n$-$\kappa X)^{\#}$.

\hspace{\parindent} In \S 2, $n$-qH-closed spaces are defined and characterized in terms of filters and covering properties, likewise with $n$-H-closed spaces. 
In \S 3, dense embeddings are constructed for $n$-Hausdorff spaces, and in \S 4 we study maximal extensions, denoted by $n$-$\kappa X$, and look at the partitions of $n$-$\sigma X\setminus X$. We use standard notation as in \cite{PW}.

\section{$n$-$q$H-closed spaces and $n$-H-closed spaces}

We recall the following.

\begin{definition}\rm
Let ${\mathcal F}$ be a filter on a space $X$. The set $\{cl_{X}(F):\;F\in{\mathcal F}\}$ is called the adherence of ${\mathcal F}$ and it is denoted by $a{\mathcal F}$. The set $\{G\in {\mathcal P}(X):\;G\supseteq F\;\hbox{for some}\;F\in{\mathcal F}\}$ is denoted by $\langle{\mathcal F}\rangle$. We say that ${\mathcal F}$ converges to $p\in X$ if ${\mathcal N}_p$ (the family of all the neighborhoods of $p$) is contained in $\langle{\mathcal F}\rangle$. Let $c({\mathcal F})$ denote the set of all convergent points of ${\mathcal F}$.

\end{definition} 

\begin{definition}\rm\cite{B}
	Let $n\in\omega$, $n\geq 2$. A space $X$ is {\bf $n$-Hausdorff} if for each distinct points $x_1,...,x_n\in X$, there exist open subsets $U_i$ of $X$ containing $x_i$ for every $i=1,...,n$ such that $x_i \in U_i$ and $\bigcap_{i=1}^{n}U_i = \varnothing$.   
\end{definition}
Clearly, we have Hausdorff when $n = 2$. 

\begin{proposition}\cite{B}	
	A space $X$ is $3$-Hausdorff iff $\Delta_3$ is closed in $Y\cup \Delta_3$, where $\Delta_3$ is the diagonal in $X^3$ and $Y=\{(x,y,z): x,y,z $ are distinct points of $X\}$.	
\end{proposition}

The following theorem gives some characterizations of $n$-Hausdorffness and also represents a generalization of the previous proposition.

\begin{theorem}\label{1}
Let $n\in\omega$, $n\geq 2$, and $X$ be a space. The following are equivalent:
	\begin{itemize}
		
		\item[(a)] $X$ is $n$-Hausdorff
		\item[(b)] if $\mathcal U$ is an open ultrafilter on $X$, then $|a\mathcal U|\leq n-1$
		\item[(c)] if $\mathcal F$ is an open filter base on $X$, $c(\mathcal F)$ contains at most $n-1$ points
		\item[(d)] $\Delta_n$ is closed in $Y\cup \Delta_n$, where $\Delta_n$ is the diagonal in $X^n$ and $Y=\{(x_1, x_2, ..., x_n): x_1, x_2, ..., x_n $ are distinct points of $X\}$.

	\end{itemize}
\end{theorem}
\proof

(a) $\Rightarrow$ (b): Assume that $|a\mathcal{U}| \geq n$, then by (a), there are open sets $\{U_a:a \in a\mathcal{U}\}$ such that $\cap_{a\in a\mathcal{U}} U_a= \varnothing$. By the ultrafilter property of $\mathcal{U}$,  $U_a \in \mathcal{U}$, and $\mathcal{U}$ has the finite intersection property (fip).  So, $\cap_{a\in a\mathcal{U}}U_a \ne \varnothing$, a contradiction.  Hence, $|a\mathcal{U}| \leq n-1$.

(b) $\Rightarrow$ (c): Let $\mathcal F$ be an open filter base on $X$. Then ${\mathcal F}$ extends to an open ultrafilter ${\mathcal U}$. By (b), $|a{\mathcal U}|\leq n-1$. Let $x\in c({\mathcal F})$. Then ${\mathcal N}_{x}\subseteq{\mathcal F}\subseteq{\mathcal U}$, so $x\in a{\mathcal U}$. Then $|c({\mathcal F})|\leq|a{\mathcal U}|\leq n-1$.

(c) $\Rightarrow$ (d): Let $x_1, x_2, ..., x_n $ be distinct points of $X$. If ${\mathcal N}_{x_1}\cup\cdots\cup {\mathcal N}_{x_n}$ has the finite intersection property, we have that 
${\mathcal N}_{x_1}\vee ....\vee {\mathcal N}_{x_n}$ is a filter converging to $x_1, x_2, ..., x_n $, contradicting (c). So there are $N_i\in {\mathcal N}_{x_i}$, where $i\in\{1,...,n\}$ such that $\bigcap _{i=1,...,n} N_i=\emptyset$. Then $N_1\times ...\times N_n$ is a neighbourhood of $(x_1, ..., x_n)$ and $(N_1\times ...\times N_n)\cap \Delta_n=\emptyset$.

(d) $\Rightarrow$ (a): Similar to (c) $\Rightarrow$ (d). 
\endproof
\begin{definition}
	Let $n\in\omega$, $n\geq 2$. An extension $Y$ of $X$ is said to be $n$-{\bf Hausdorff except for $X$}  if for points $p \in Y\backslash X$ and $q_1,... , q_{n-1} \in Y$, there are open sets $U, V_1, ..., V_{n-1}$ in $Y$ such that $p \in U$ and $q_i \in V_i$, $i=1,..., n-1$ and $U\cap V_1\cap ,,, \cap V_{n-1}=\emptyset$.
\end{definition}

\begin{theorem}\label{new}
Let $n\in\omega$, $n\geq 2$. The following are equivalent for any space $X$:
	\begin{enumerate}
				\item[(a)] For every open filter ${\mathcal F}$ on $X$, $|a{\mathcal F}|\geq n\hbox{-}1$;
		\item[(b)] For every $A\in[X]^{<n\hbox{-}1}$ and for family of open subsets ${\mathcal U}$ of $X$ such that $X\setminus A\subseteq \bigcup{\mathcal U}$, there exists ${\mathcal V}\in[{\mathcal U}]^{<\omega}$ such that $X=\bigcup_{V\in{\mathcal V}}\overline{V}$.
	\item [(c)] $X$ is closed in every extension of $X$ that is $n$-Hausdorff expect for $ X$
	\end{enumerate}
\end{theorem}

\proof
(a) $\Rightarrow$ (b) Let $A\in[X]^{< n\hbox{-}1}$. Suppose there exists $\mathcal{C}$ a family of open subsets of $X$ containing $X\setminus A$ such that for every finite subfamily ${\mathcal V}$ of ${\mathcal C}$, $X\neq \bigcup_{V\in{\mathcal V}}\overline{V}$. $F=\{U:\;U\;\hbox{is open and}\;U\supseteq X\setminus\bigcup_{V\in{\mathcal V}}\overline{V}\;\hbox{for some finite subfamily}\;{\mathcal V}\;\hbox{of}\;{\mathcal C}\}$ is an open filter. As $a{\mathcal F}=\bigcap_{U\in{\mathcal F}}\overline{U}\subseteq \bigcap_{{\mathcal V}\in[{\mathcal C}]^{<\omega}}\overline{X\setminus \bigcup_{V\in{\mathcal V}}\overline{V}}\subseteq \bigcap_{V\in{\mathcal V}}\overline{X\setminus\overline{V}}\subseteq X\setminus {\mathcal C}=X\setminus (X\setminus A)=A$, we have $|a{\mathcal F}|<n\hbox{-}1$. A contradiction.

(b) $\Rightarrow$ (a) Let ${\mathcal F}$ be an open filter on $X$ such that $|a{\mathcal F}|<n\hbox{-}1$. ${\mathcal W}=\{X\setminus\overline{U}:\;U\in{\mathcal F}\}$ is a family open sets covering $X\setminus a{\mathcal F}$, by (b) there exists a finite subset $A$ of ${\mathcal F}$ such that $X=\overline{\bigcup_{U\in A}X\setminus\overline{U}}=\bigcup_{U\in A}X\setminus int(\overline{U})=X\setminus\bigcap_{U\in A} int(\overline{U})\subseteq X\setminus \bigcap_{U\in A}U$, then $\bigcap_{U\in A}U=\emptyset$. A contradiction.

(a) $\Rightarrow$ (c) Assume every filter ${\mathcal F}$ has $|a{\mathcal F}|\geq n-1$. By the way of contradiction, let $Y$ be an extension of $X$ such that $X\neq cl_{Y}(X)$ and $Y$ is $n$-Hausdorff except for $X$. Let $p\in cl_{Y}(X)\setminus X$. Then ${\mathcal F}=\{U\cap X:\;p\in U\in\tau(Y)\}$ is an open filter on $X$. Extend ${\mathcal F}$ to an open ultrafilter ${\mathcal U}$. Then $|c{\mathcal U}|\geq n-1$. Let $x_{1},...,x_{n-1}\in c\mathcal U$. As $Y$ is $n$-Hausdorff except for $X$, there exists open sets in $Y$ $U_{i}$, $U$ such that $x_{i}\in U_{i}$, $i=1,...,n-1$, $p\in U$ and $\bigcap_{i=1,...,n}^{n-1}U_{i}\cap U=\emptyset$. We can observe that $X\cap U_{i}\in{\mathcal U}$ for every $i=1,...,n-1$ and $U\cap X\in{\mathcal U}$, then $\bigcap_{i=1,...,n}^{n-1}(X\cap U_{i})\cap (U\cap X)\in{\mathcal U}$. Therefore $\bigcap_{i=1,...,n}^{n-1}U_{i}\cap U\neq\emptyset$; a contradiction. 

(c) $\Rightarrow$ (a) Let $\mathcal F$ be an open filter such that $|a{\mathcal F}|<n-1$, consider the extension $Y=X\cup \{{\mathcal F}\}$, and define a set $U\subseteq Y$ to be open if $U\cap X$ is open in $X$ and ${\mathcal F}\in U$ implies $U\cap X\in{\mathcal F}$. It is easy to prove that $Y$ is Hausdorff (hence $n$-Hausdorff) except for $X$ but $X$ is not closed in $Y$ (see also  Proposition 4.8 (b) in \cite{PW}).

\begin{definition}
Let $n\in\omega$, $n\geq 2$.	A space is $n$-qH-closed if $X$ satisfies one (and hence all) of the conditions of Theorem \ref{new}.
\end{definition}

\begin{definition}\rm
$2$-qH-closed spaces are called qH-closed spaces.
\end{definition}

Clearly we have:
\begin{proposition}\rm
Let $n\in\omega$, $n\geq 2$. Every $(n+1)$-qH-closed space is $n$-qH-closed.
\end{proposition}

We can notice the following behaviour of $n$-qH-closedness with respect to regular closed subsets.

\begin{proposition}Let $n\in\omega$, $n\geq 2$. If $X$ is a $n$-qH-closed space and $U$ is an open subset of $X$ then $cl_{X}(U)$ is $n$-qH-closed.
\end{proposition}
\proof We put $A=cl_{X}(U)$ and let ${\mathcal F}$ be an open filter of $A$. As ${\mathcal G}=\{W\subseteq X:\;W\;\hbox{is an open subset of}\;X,\;W\supseteq F\cap U,\;\forall F\in{\mathcal F}\}$ is an open filter on $X$, we have $|a_X{\mathcal G}|\geq n-1$. We have $a_X{\mathcal G}\subseteq a_A{\mathcal F}$ then $|a_A{\mathcal F}|\geq n-1$. By Theorem \ref{new}, $A$ is $n$-qH-closed.
\endproof

The following represents a generalization of the concept of H-closedness. 

\begin{definition}\rm
Let $n\in\omega$, $n\geq 2$.	An n-Hausdorff space $X$ is called {\bf n-H-closed} if  $X$ is closed in every n-Hausdorff space $Y$ in which $X$ is embedded.  
\end{definition}

Clearly a space $X$ is 2-H-closed iff $X$ is H-closed. 
  
Obviously, for every $n\in\omega$, $n\geq2$, every $n$-Hausdorff space is $n+1$-Hausdorff. In Example 4 in \cite{B}, a $(n+1)$-Hausdorff non $n$-Hausdorff space is given, for every $n\in\omega$. Also we have the following

\begin{example}\label{esempio1} An H-closed non 3-H-closed space.
 \end{example} 
  Let $Y=3\cup (\omega \times 3)$, where $3=\{0,1,2\}$, be the space topologized as follows: all points from $\omega\times 3$ are isolated; a basic neighborhood of $i\in 3$ takes the form $U(i,N)=\{i\}\cup\{(m,j): j\neq i, m\geq N\}$, where $N\in \omega$. $Y$ is $3$-Hausdorff non Hausdorff. Consider the subspace $X=2\cup (\omega \times 2)\subset Y$.
  $X$ is compact and Hausdorff, then H-closed but $X$ is not closed in $Y$.

 The following is easy to verify.
 
\begin{proposition}
An n-Hausdorff space  $X$ is n-H-closed iff $X$ is closed in every n-Hausdorff space $Y$ in which $X$ is embedded where $|Y\backslash X| =1.$
\end{proposition}

Now we consider some characterizations of $n$-H-closed spaces.

\begin{theorem}\label{2} Let $n\in\omega$, $n\geq 2$. For a n-Hausdorff space $X$ the following are equivalent
	\begin{enumerate}
		\item [(a)] $X$ is n-H-closed 
		\item [(b)] for each open ultrafilter $\mathcal{U}$ on $X$, $|a\mathcal{U}| = n-1$
		\item [(c)] $X$ is $n$-qH-closed
		\item [(d)] for every ${\mathcal F}_{1},\cdots,{\mathcal F}_{n-1}$ open filters on $X$ such that $\bigcup_{1=1}^{n-1}\mathcal F_{i}$ has the finite intersection property, we have $|a({\mathcal F}_{1}\vee...\vee{\mathcal F}_{n-1})| \geq n-1$
\end{enumerate}
\end{theorem}

\proof 
(b) $\Rightarrow$ (a) Suppose for each open ultrafilter $\mathcal{U}$ on $X$, $|a\mathcal{U}| = n-1,$ and   $X$ is embedded in an n-Hausdorff space $Y$ such that $Y\backslash X = \{p\}$. Assume $p \in cl_YX$.  Let $\mathcal{B}_p = \{U \cap X: p \in U \in \tau(Y)\}$ and $\mathcal{U}_p$ be an open ultrafilter in $X$ containing $\mathcal{B}_p$.  By hypothesis, $|a_X(\mathcal{U}_p)| = n-1$.  Let $\mathcal{W}_p$ be the unique open ultrafilter in $cl_Y(X)$ containing $\mathcal{U}_p$ and $T \in \mathcal{W}_p$.  If $p \in V \in \tau(Y)$, then $V\cap X \in \mathcal{B}_p$ and $(V \cap X) \cap (T\cap X) \ne \varnothing$ for each $T \in \mathcal{W}_p$.  Thus, $V \cap T \ne \varnothing$. This shows that $p \in a(\mathcal{W}_p)$ and $|a(\mathcal{W}_p)| = n$, a contradiction.   

(a) $\Rightarrow$ (b) Suppose $X$ is n-H-closed and assume there is an open ultrafilter $\mathcal{U}$ on $X$ such that  $|a\mathcal{U}| < n-1$. Let $Y = X \cup \{p\}$ where $p \not\in X$.  We define the topology on $Y$ by $U$ is open in $Y$ if $U \cap X$ is open in $X$ and if $p \in U$, then $U\cap X \in \mathcal{U}$.  Thus, $p \in cl_Y(X)$.  To obtain a contradiction, we will show that $Y$ is n-Hausdorff.  Let $A \subseteq Y$ such that $|A|= n$.  If $A \subseteq X$, there are open sets $U_a$ in $X$ for each $a \in A$ such that $\bigcap_{a \in A}A_a = \varnothing$ and as $X$ is open in $Y$, $A_a$ is also open in $Y$. Otherwise, $p \in A$ and $|A\backslash\{p\}| = n-1$.  As $|a\mathcal{U}| < n-1$, there is $q \in  A\backslash a\mathcal{U}$ and $q \ne p$.  There is open set $V$ in $X$ and $U \in \mathcal{U}$ such that $q \in V$ and $V \cap U = \varnothing$.  Both $V$ and $\{p\} \cup U$ are disjoint open sets in $Y$ and $\{p,q\} \subseteq A$.  This completes the proof that $Y$ is n-Hausdorff, a contradiction.  

(a) $\Rightarrow$ (c) Let ${\mathcal F}$ be an open filter on $X$. As ${\mathcal F}$ is contained in some open ultrafilter ${\mathcal U}$ on $X$, we have $a{\mathcal F}\supseteq a{\mathcal U}$ and by (a) we have $|a{\mathcal F}|\geq n\hbox{-}1$.

(c) $\Rightarrow$ (a) Let ${\mathcal U}$ be an open ultrafilter on $X$, by (b) we have $|{\mathcal U}|\geq n\hbox{-}1$. By Theorem \ref{1} $|{\mathcal U}|\leq n\hbox{-}1$ and $|{\mathcal U}|= n\hbox{-}1$.

(b) $\Rightarrow$ (d) The open filter ${\mathcal F}_{1}\vee...\vee{\mathcal F}_{n-1}$ is contained in an open ultrafilter $\mathcal{U}$ and $|a\mathcal{U}| = n-1$ by (b).  As $|a({\mathcal F}_{1}\vee...\vee{\mathcal F}_{n-1})| \geq|  a\mathcal{U}|$, it follows that $|a({\mathcal F}_{1}\vee...\vee{\mathcal F}_{n-1}) |\geq n-1$. 

(d) $\Rightarrow$ (c) is obvious.
\endproof

\hspace{\parindent} Note that the space $X$ of Example \ref{esempio1} is an H-closed not 3-H-closed space. We can also observe that an H-closed space $X$ is {\bf not} $n$-H-closed for $n >2$. As every open ultrafilter on an $n$-H-closed space $X$ has $n-1 $ adherence points but in an H-closed space, every open ultrafilter has an unique adherent point.

\begin{comment} Let $n\in\omega$, $n\geq 2$. For a n-Hausdorff space $X$, an open ultrafilter $\mathcal{U}$ on $X$ is said to be {\bf full} if $|a\mathcal{U}| = n-1$.  So, by Theorem \ref{2}, a n-Hausdorff space $X$ is n-H-closed iff every open ultrafilter on $X$ is full. Consider the simple space $X = \omega \cup\{a,b\}$ where a set $U \subseteq X$ is defined to be open   if $a \in U$ or $b \in U$, then $U\cap \omega$ is cofinite in $\omega$, and each point $n\in\omega$ is isolated. The compact space $X$ is 3-Hausdorff but not 3-H-closed for if $\mathcal{U}$ is an open ultrafilter containing the open set $\{1\}$, then $a\mathcal{U} = \{1\}$ and $\mathcal{U}$ is not full. 
\end{comment}

\section{Embeddings}
\hspace{\parindent} In 1924, Alexandroff and Urysohn asked if every Hausdorff space can be embedded densely in an H-closed space. Kat\v{e}tov and Stone answered this (see page 307 in \cite{PW}). We ask a similar question: can every n-Hausdorff space be densely embedded in an n-H-closed space?
We will use open ultrafilters in answering this question in the affirmative. We start the construction of an extension of $X$ by adding one open ultrafilter that is not full.  This is a modification of the proof of (a) $\Longleftrightarrow$ (b) in Theorem \ref{2}.

Recall the following definitions and constructions given in \cite{PW}.

\hspace{\parindent} Let $Y$ be an extension of a space $X$.  For $p \in Y$, let $O^p = \{U \cap X: p \in U \in \tau(Y)\}$ and for $U \in \tau(X)$, let $oU = \{p \in Y: U \in O^p\}$.  Note that for $U,V \in \tau(X)$, $o(U\cap V) = oU \cap oV$, $o(\varnothing) = \varnothing$, and $oX =   Y$. So, $\{oU:U \in \tau(X)\}$ forms a basis for a topology, denoted as $\tau^{\#}(Y)$, on $Y$.  Denote by $\tau^+(Y)$ the topology on $Y$ generated by the base  ${\mathcal B}=\{U\cap \{p\}: U\in O^p$ and $p\in Y\}$.  We have $\tau^{\#}(Y) \subseteq \tau(Y) \subseteq \tau^+(Y)$, and $Y^\#$(resp. $Y^+$) is used to denote the set $Y$ with $\tau^{\#}(Y)$ (resp. $\tau^+(Y)$). $Y^+$ is called {\bf simple extension} of the space $X$ and $Y^{\#}$ is called {\bf strict extension} of the space $X$.

\begin{proposition}\label{B}
	Let $Y$ be an extension of a space $X$. 
	Then $Y^+$ and $Y^{\#}$ are extensions of $X$ with the following properties:
	\begin{enumerate}
		\item [(a)]$(Y^+)^+ = Y^+$ and $(Y^{\#})^{\#} = Y^{\#}$;
		\item [(b)]$Y$ is qH-closed iff $Y^+$ is qH-closed iff $Y^{\#}$ is qH-closed;
		\item [(c)]$Y$ is Hausdorff iff $Y^+$ is Hausdorff iff $Y^{\#}$ is Hausdorff;
		\item [(d)]If $U \in \tau(X)$, $cl_YU = cl_{Y^+}U = cl_{Y^{\#}}U$;
		\item [(e)]If $p \in Y$, $O^p_Y = O^p_{Y^+} = O^p_{Y^{\#}}$;
		\item [(f)]If $\sigma$ is a topology on the set $Y$, $(Y,\sigma)$ is an extension of $X$, and $O^p_{\sigma} = O^p_Y$ for all $p \in Y$, then $\tau^{\#}(Y) \subseteq \sigma \subseteq \tau^+(Y)$;
		\item [(g)]The subspace $Y^+\backslash X$ is discrete.
	\end{enumerate}
\end{proposition}

\begin{proposition}\label{C}
	Let $Y$ be an extension of a space $X$, $\mathcal{U}$ an open ultrafilter on $X$, and $\mathcal{U}^e =  \{V \in \tau(Y): V \cap X \in \mathcal{U}\}$. Then: \begin{enumerate}
		\item [(a)] $\mathcal{U}^e$ is an open ultrafilter on $Y$,
		\item [(b)]if $\mathcal{V}$ is an open ultrafilter on $Y$ such that $\{V\cap X:V \in \mathcal{V}\} \subseteq \mathcal{U}$, then $\mathcal{U}^e = \mathcal{V}$ and $a\mathcal{V} = \bigcap_{\mathcal{V}}cl_YU \supseteq a\mathcal{U}$.
	\end{enumerate}
\end{proposition}

\begin{proposition}\label{0.5}
Let $n\in\omega$, $n\geq 2$, and $Y$ be a $n$-Hausdorff extension of a space $X$.  Then both the simple extension $Y^+$ and strict extension $Y^{\#}$ of $X$ are $n$-Hausdorff as well as $X$.
\end{proposition}
\proof
As $\tau(Y) \subseteq \tau(Y^+)$, it is immediate that $Y^+$ is also $n$-Hausdorff. Likewise, it is easy to verify that $X$ is also $n$-Hausdorff. To show that 
$Y^{\#}$ is $n$-Hausdorff, let $A \subseteq Y^{\#}$ such that $|A| = n$.  As $Y$ is $n$-Hausdorff extension, there is a family of open sets $\{V_a:a \in A\}$ such that 
$a \in V_a$ and $\bigcap_{a \in A}V_a=\varnothing$. So, $\bigcap_{a \in A}(V_a\cap X)=\varnothing$.  Thus, $\varnothing =o(\varnothing) = o(\bigcap_{a \in A}(V_a\cap X)) = \bigcap_{a \in A}o(V_a\cap X) $.  This shows that $Y^{\#}$ is $n$-Hausdorff.
\endproof

\begin{proposition}\label{cor1}
Let $n\in\omega$, $n\geq 2$.	If $Y$ is an extension of $X$, then $Y$ is $n$-H-closed iff $Y^{\#}$ is $n$-H-closed iff $Y^+$ is $n$-H-closed.
\end{proposition}
\proof
This is an application of Proposition \ref{0.5}, Proposition \ref{C}, and Proposition \ref{B}(d).
\endproof

\begin{ECT}\rm Let $n\in\omega$, $n\geq 2$, and $X$ be an n-Hausdorff space and $\mathcal{U}$ an open ultrafilter on $X$ such that $|a{\mathcal{U}}| < n-1$.  Let $k = n-1-|a{\mathcal{U}}|$, and  $Y = X \cup\{p_1,p_2, \cdots, p_k\}$ where $\{p_1,p_2, \cdots, p_k\} \cap X = \varnothing$. A set $V$ is defined to be open in $Y$ if $V \cap X$ is open in $X$ and if $p_i \in V$ for $1 \leq i \leq k$, $V \cap X \in \mathcal{U}$.  Now, $Y \backslash X = \{p_1,p_2, \cdots, p_k\}$, and a basic open set containing $p_i$ is $\{p_i\} \cup T$ where $T \in \mathcal{U}$. It is straightforward to verify the following result using the technique developed in the proof of Theorem \ref{2}.
\end{ECT}

\begin{lemma} Let $n\in\omega$, $n\geq 2$, $X$ be an n-Hausdorff space and $Y$ the space defined in the above construction. The space $Y$ is an n-Hausdorff space that contains $X$ as a dense subspace and if $\mathcal{V}$ is the unique open ultrafilter on $Y$ containing $\mathcal{U}$, $|a_Y(\mathcal{V})| = n-1$.  $\mathcal{W}$ is an open ultrafilter on $Y$ iff $\mathcal{W} = \mathcal{V}$ or $\mathcal{W}$ is an open ultrafilter on $X$ other than $\mathcal{U}$.
\end{lemma}

\begin{ECT}\rm Let $n\in\omega$, $n\geq 2$, $X$ be an n-Hausdorff space and  $\mathfrak{U} = \{\mathcal{U}:\mathcal{U} \text{ is an 
open}$ ultrafilter such that  $|a{\mathcal{U}}| < n-1\}$.  That is, $\mathfrak{U}$ is the set of open ultrafilters on $X$ that are not full.  We indexed $\mathfrak{U}$ by $\mathfrak{U} = \{ \mathcal{U}_{\alpha}: \alpha \in |\mathfrak{U}|\}$.  For each $\alpha \in |\mathfrak{U}|$, let $k\alpha = n - 1 - |a\mathcal{U}_{\alpha}|$ and $\{p_{\alpha i}: 1 \leq i \leq k\alpha\}$ a set of distinct points disjoint from $X$.  Let $Y  = X \cup \{p_{\alpha i}:1 \leq i \leq k\alpha,\; \alpha \in |\mathfrak{U}| \}$. A set $V$ is defined to be open in $Y$ if $V \cap X$ is open in $X$ and if $p_{\alpha i} \in V$ for $1 \leq i \leq k\alpha$, $V \cap X \in \mathcal{U}_{\alpha}$.  The space $Y$ is an extension of $X$ such that $Y$ is an n-Hausdorff space in which every open ultrafilter on $Y$ is full.  That is, $Y$ is n-H-closed.  These properties are summarized in the following result.
\end{ECT}

\begin{theorem} Let $n\in\omega$, $n\geq 2$. An n-Hausdorff space can be densely embedded in an n-H-closed space.
\end{theorem}

In the next result, we  present the basic extension properties of $O^p$ and $oU$ and $cl_YU$ for $U \in \tau(X)$ for this particular extension $Y$.

\begin{proposition}\label{1.5}
Let $n\in\omega$, $n\geq 2$, $X$ be $n$-Hausdorff space and $Y$ the construction of an $n$-H-closed simple extension of $X$.  Then:
\begin{enumerate}

\item [(a)]For $p = p_{\alpha i} \in Y\backslash X$, $O^{p_{\alpha i}} = \{V \cap X: p \in V \in \tau(Y)\} = \mathcal{U}_{\alpha}$;
\item [(b)]For $p \in X$, $O^{p} = \{V \cap X: p \in V \in \tau(Y)\} =\{U \in \tau(X): p \in U\}$;
\item [(c)]For $U \in \tau(X)$, $oU = \{p \in Y:U \in O^p\} = U \cup \{p_{\alpha i}: U \in \mathcal{U}_{\alpha}, 1\leq i\leq k\alpha \}$;
\item [(d)] For $U \in \tau(X)$, $cl_YU = cl_XU \cup oU$. \end{enumerate}
\end{proposition}
\proof
The proof of (a), (b), and (c) are straightforward and left to the reader.  To show (d), let $p_{\alpha i} \in cl_YU\backslash X$.  Then $V \cap U \ne \varnothing$ for each $V \in \mathcal{U}_{\alpha}$ implying $U \in \mathcal{U_{\alpha}}$ and that $p_{\alpha i} \in oU$ by (a) and (c).  Note that $cl_YU\cap X = cl_XU$.  Thus, $cl_YU \subseteq  cl_XU \cup oU$.  Clearly, $cl_XU \subseteq cl_YU$ and by (a) and (c), $oU \subseteq cl_YU$.  This completes the proof of (d).
\endproof

\begin{comment}\label{10}\rm Let $n\in\omega$, $n\geq 2$. The above construction of an n-H-closed extension works whenever $X$ is n-Hausdorff. Consider the example of a 3-Hausdorff space  
$X = \omega \cup\{a,b\}$  described in Remark 1 that is not 3-H-closed. Notice that every open ultrafilter on $X$ containing only infinite subsets is an open ultrafilter $\mathcal{U}$ on $\omega$ and $a\mathcal{U} = \{a,b\}$. The only open ultrafilters remaining are those containing a single point set like $\{n\}$. 
To construct a 3-H-closed space $Y$ that contains $X$ as a dense subspace let $Y = X \cup \{c_n: n \in \omega\}$.  A set $U \subseteq Y$ is defined to be open if $U \cap X$ is open in $X$ and if  $c_n \in U$, then  $n \in U$.   
The space $Y$ is 3-Hausdorff, and $X$ is dense in $Y$.  Let $\mathcal{U}$ be an open ultrafilter on $X$ and $\mathcal{V} = \{U  : U$ is open in Y and $U \cap X \in \mathcal{U}\}$ the unique open ultrafilter on $Y$ that contains $\mathcal{U}$.  If $U$ is infinite for each $U \in \mathcal{U}$, then $\{a,b\} = a\mathcal{V}$.  If some $U$ is finite for $U \in \mathcal{U}$, then there is an unique $n \in \omega$ such that $U = \{n\}$, and it follows that $\{c_n, n\} = a\mathcal{V}$. Thus, $Y$ is a 3-H-closed space. The space $Y$ provides an example that shows the size of $a\mathcal{F}$ can be quite large.
For $n \in \omega$, the open ultrafilter $\mathcal{U}_n = \{U \in \tau(Y):n \in U\}$ has the property that $|a\mathcal{U}_n| = 2$.  However, the open filter $\mathcal{F} = \bigcap_{1\leq i \leq n}\mathcal{U}_n$ on $Y$ has the property that $|a\mathcal{F}| = 2n$.  This shows that the adherence of an open filter $\mathcal{F}$ in an $n$-H-closed space can be large.

\end{comment}

\section{Theory of $n$-H-closed Extensions.}  

\begin{theorem}\label{7}
Let $n\in\omega$, $n\geq 2$, $X$ be n-Hausdorff space and $Y$ the simple, n-H-closed extension of $X$ constructed in Extension Construction Technique 2.  If $Z$ is an n-H-closed extension of $X$, there is a continuous surjection $f:Y \rightarrow Z$ such that $f(x) = x$ for all $x \in X$.
\end{theorem}

\proof
Let $\hat{u} = \{\mathcal{U}:$ $\mathcal{U}$ is an open ultrafilter on $X$ such that $|a\mathcal{U}| < n-1\}$.  For $\mathcal{U} \in \hat{u}$, the collection  $\mathcal{V}_{\mathcal{U}} =\{V \in \tau(Z): V\cap X \in \mathcal{U}\}$ is an open ultrafilter on $Z$. Note that $a_Z\mathcal{V}_{\mathcal{U}}\backslash X = a_Z\mathcal{V}_{\mathcal{U}}\backslash a_X\mathcal{U} = \{p_1,p_2,\cdots,p_k\}$ where $k = n-1-|a_X\mathcal{U}|$ and $Z = X \cup \bigcup_{\mathcal{U} \in \hat{u}}a_Z\mathcal{V}_{\mathcal{U}}\backslash a_X\mathcal{U} .$\\
For $\mathcal{U} \in \hat{u}$, let $\mathcal{V'}_{\mathcal{U}} =\{V' \in \tau(Y): V'\cap X \in \mathcal{U}\}$ is an open ultrafilter on on $Y$. The set $a_Y\mathcal{V'}_{\mathcal{U}}\backslash a_X\mathcal{U} = \{q_1,q_2,\cdots,q_k\}$ where $k = n-1-|a_X\mathcal{U}|$. Note that $Y = X \cup \bigcup_{\mathcal{U} \in \hat{u}}a_Y\mathcal{V'}_{\mathcal{U}}\backslash a_X\mathcal{U}.$\\
Define  $f_{\mathcal{U}}(q_i) = p_i$ for $1 \leq i \leq k$; $f_{\mathcal{U}}$ is a bijection.  Define $f:Y\rightarrow Z$ as follows: for $x \in X$, $f(x) = x$ and for $q_i \in \mathcal{V'}_{\mathcal{U}}$, $f(q_i) = f_{\mathcal{U}}(q_i) = p_i$. The function $f$ is onto, $f(x) = x$ for $x \in X$, but not necessarily one-to-one. As $X$ is open in $Y$, $f$ is continuous for $x \in X$.  For $\mathcal{U} \in \hat{u}$ and $q_i \in a_Y\mathcal{V'}_{\mathcal{U}}\backslash a_X\mathcal{U}$, $f(q_i) = f_{\mathcal{U}}(q_i) = p_i$; let $p_i \in V$ for some $V$ open in $Z$.  Then $V \cap X \in \mathcal{U}$ and   $f[\{q_i\} \cup (V \cap X)] \subseteq V$.  As  $\{q_i\} \cup (V \cap X)$ is open in $Y$, it follows that  $f$ is continuous at $q_i$. \qed
 
\begin{comment}\rm
Let $n\in\omega$, $n\geq 2$, $S$ and $T$ be $n$-H-closed extensions of a $n$-Hausdorff space $X$.  We say $S$ is {\bf projectively larger} than $T$ if there is a continuous surjection $f:S\rightarrow T$ such that $f(x) = x$ for $x \in X$. This projectively larger function may not be unique.
\end{comment}
The proof of Theorem \ref{7} shows that the $n$-H-closed extension $Y$ of $X$ is projective larger than every $n$-H-closed extension of $X$. The space $Y$ from Theorem \ref{7} has an interesting uniqueness property as noted in the next result.

\begin{theorem}\label{8}
Let $n\in\omega$, $n\geq 2$, $X$ be an $n$-Hausdorff space and $Y$ the $n$-H-closed extension of $X$ described in Theorem \ref{7}.  Let $f:Y\rightarrow Y$ be a continuous surjection such that $f(x) = x$ for all $x \in X$.  Then $f$ is a homeomorphism.
\end{theorem}
 
\proof From Theorem \ref{7}, we have that  $Y = X \cup \bigcup_{\mathcal{U} \in \hat{u}}a_Y\mathcal{V'}_{\mathcal{U}}\backslash a_X\mathcal{U}$ where if $\mathcal{S}$, $\mathcal{T} \in \hat{u}$ are distinct open ultrafilters on $X$, then 
$ (a_Y\mathcal{V'}_{\mathcal{S}}\backslash a_X\mathcal{S}) \cap (a_Y\mathcal{V'}_{\mathcal{T}}\backslash a_X\mathcal{T}) = \varnothing$.  Let $\mathcal{U} \in \overline{u}$ and $\mathcal{V}$ be the unique open ultrafilter on $Y$ such that $\mathcal{U} \subseteq \mathcal{V}$.  Then $f[a\mathcal{V}] = f[\bigcap_{V \in \mathcal{V}}cl_YV] \subseteq \bigcap_{V \in \mathcal{V}}f[cl_YV] =\bigcap_{V \in \mathcal{V}}f[cl_Y(V\cap X)] \subseteq \bigcap_{V \in \mathcal{V}}cl_Yf[V\cap X] = \bigcap_{V \in \mathcal{V}}cl_YV = a\mathcal{V}$.  Now, $f[a\mathcal{V}] = f[a\mathcal{V}\backslash a\mathcal{U} \cup a\mathcal{U}] =f[a\mathcal{V}\backslash a\mathcal{U}] \cup a\mathcal{U} \subseteq a\mathcal{V}\backslash a\mathcal{U} \cup a\mathcal{U} $.  If  $\mathcal{T} \in \hat{u}$, $\mathcal{T} \ne \mathcal{U}$, and $\mathcal{V'}$ is the unique open ultrafilter on $Y$ such that $\mathcal{T} \subseteq \mathcal{V'}$, then $f[a\mathcal{V'}] \subseteq a\mathcal{V'}\backslash a\mathcal{T} \cup a\mathcal{T}$.  But $a\mathcal{V}\backslash a\mathcal{U}  \cap  a\mathcal{V'}\backslash a\mathcal{T} = \varnothing.$. Thus, if $f$ is onto, $f[a\mathcal{V}\backslash a\mathcal{U}] = a\mathcal{V}\backslash a\mathcal{U}$.  Since $a\mathcal{V}\backslash a\mathcal{U}$ is finite, then $f|_{a\mathcal{V}\backslash a\mathcal{U}}:a\mathcal{V}\backslash a\mathcal{U}\rightarrow a\mathcal{V}\backslash a\mathcal{U}$ is also one-to-one. This completes the proof that $f$ is a bijection.  To show that $f$ is open first note that for an open set $V \subseteq X$, $f[V] = V$ is open. Let $p \in a\mathcal{V}\backslash a\mathcal{U}$ and $p \in V \in \mathcal{V}$. 
Then $ \{p\} \cup V\cap X$ is a basic open set in $Y$ containing $p$. Now, $f[\{p\} \cup V\cap X] = \{f(p)\} \cup f[V \cap X] = \{f(p)\} \cup (V \cap X)$ is a basic open set in $Y$ containing $f(p)$.   This completes the proof that $f$ is open. \qed

\begin{comment}\rm
In the setting of Hausdorff spaces, projective maximums and projectively larger functions (defined in \cite{PW}) are unique.   Sometimes this is  a problem in non-Hausdorff spaces. From Theorem \ref{8}, we see that a form of uniqueness for the class of $n$-H-closed extensions is possible. We extend the definition of projective maximum in \cite{PW} as follows:  Let $\mathcal{E}(X)$ be a collection of extensions of a space $X$. We say $Y \in \mathcal{E}(X)$ is a {\bf projective maximum} of $X$ if $Y$ is projective larger than each $Z \in \mathcal{E}(X)$ and if $f:Y \rightarrow Y$ is a continuous surjection such that $f(x) = x$ for each $x \in X$, then $f$ is a homeomorphism.  By Theorem \ref{8}, the $n$-H-closed space $Y$ constructed in Theorem \ref{7} for a $n$-Hausdorff space $X$ is a projective maximum.  We denote this projective maximum by $n$-$\kappa X$ and call it the {\bf $n$-Kat\v etov extension} of $X$.
\end{comment}

\subsection{Fomin H-closed Extension} Let $n\in\omega$, $n\geq 2$, $X$ be a $n$-Hausdorff space and $Y$ denote $n$-$\kappa X$ from Theorem \ref{7}.  In the setting of Hausdorff spaces, the combination of the Kat\v etov and Fomin extensions provide the major support  for the theory of H-closed extensions. In the class of $n$-Hausdorff spaces,  the {\bf Fomin} extension is defined as $Y^\#$ and denoted as $n$-$\sigma X$.  By Proposition \ref{cor1}, $n$-$\sigma X$ is $n$-H-closed, and by Proposition \ref{B}(f), the identity function $id:n$-$\kappa X \rightarrow n$-$\sigma X$ is a continuous bijection such that $id(x) = x$ for all $x \in X$.  Recall that for spaces $S,T$, a function $f:S \rightarrow T$ is $\theta$-continuous if for each $p \in S$ and open set $V\in \tau(T)$ such that $f(p) \in V$, there is an open set $U \in \tau(S)$ such that $f[cl_SU] \subseteq cl_TV$.  By Proposition \ref{B}(d), it follows that $id:n$-$\sigma X \rightarrow n$-$\kappa X$  is $\theta$-continuous.  If $Z$ is an $n$-H-closed extension of $X$, there is a continuous surjection $f:n$-$\kappa X \rightarrow Z$ such that $f(x) = x$ for $x \in X$.  So, the composition $f\circ id:n$-$\sigma X \rightarrow Z$ is $\theta$-continuous, onto, and $(f\circ id)^{\leftarrow}[Z\backslash X] = n$-$\sigma X\backslash X$.  The next result presents some interesting properties of $n$-$\sigma X\backslash X$.

\begin{theorem}\label{9}
Let $n\in\omega$, $n\geq 2$, $X$ be an $n$-Hausdorff space and $Y^{\#}$ denote the Fomin extension $n$-$\sigma X$ of $X$.  Then: \newline
{\phantom -}\hskip 6.5mm(a) For $U \in \tau(X)$, $cl_{Y^{\#}}oU = cl_{Y^{\#}}U = cl_XU \cup oU$ and $cl_{Y^{\#}}U \backslash X =  oU\backslash X$.  \newline
{\phantom -}\hskip 6.5mm(b) For $V \in \tau(Y^{\#}), $ $cl_{Y^{\#}}V = cl_{Y^{\#}}(V \cap X) = cl_X(V\cap X) \cup o(V\cap X)$\newline
{\phantom -}\hskip 6.5mm(c) The subspace $n$-$\sigma X\backslash X$ has a clopen basis (and is zero-dimensional).\newline
{\phantom -}\hskip 6.5mm(d) If $K \subseteq n$-$\sigma X\backslash X$ is closed in $n$-$\sigma X$, then $K$ is compact. \newline
{\phantom -}\hskip 6.5mm(e) If $A \subseteq X$ is closed and nowhere dense in $X$, then $A $ is closed in  $n$-$\sigma X$.
\end{theorem}
\proof
For (a), let $U \in \tau(X)$.  Then $cl_{Y^{\#}}U = cl_{Y}U =cl_XU \cup oU$ by Propositions \ref{B}(d) and \ref{1.5}(d). (b) and (c) follow from (a). To show (d), let $\{o(U_a): a \in A\}$ be an open cover of $K$.  As $K$ is closed in $n$-$\sigma X$, there is a family of basic open sets $\{o(U_b): b \in B\}$ such that $o(U_b) \cap K = \varnothing$ for each $b \in B$.  As $n$-$\sigma X$ is $n$-H-closed by Proposition \ref{cor1}, then $n$-$\sigma X$ is qH-closed by Theorem \ref{2}.  There are finite families $A' \subseteq A$ and $B' \subseteq B$ such that  $n$-$\sigma X = \bigcup_{a\in A'}cl_{Y^{\#}}o(U_a) \cup \bigcup_{b\in b'}cl_{Y^{\#}}o(U_b) = \bigcup_{a\in A'}cl_X U_a \cup o(U_a) \cup \bigcup_{b\in b'}cl_X U_b \cup o(U_b)$.  As $K \subseteq  n$-$\sigma X\backslash X$, we have that $K \subseteq \bigcup_{a\in A'}cl_X o(U_a)$. To show (e), first note that $X\backslash A$ is open and dense in $X$.  So, $o(X\backslash A) = (X\backslash A) \cup n$-$\sigma X\backslash X$, and $n-\sigma X\backslash o(X\backslash A) = n-\sigma X\setminus (X\setminus A)\cap n-\sigma X\setminus (n-\sigma X\setminus X)= n-\sigma X\setminus (X\setminus A)\cap X=A $.   \qed

\begin{comment}\rm
Let $n\in\omega$, $n\geq 2$, and $X$ be an $n$-Hausdorff space.  By Theorem 9(c), the space $n$-$\sigma X\backslash X$ is completely regular (Tychonoff without being Hausdorff) in addition to being zero-dimensional. In fact, by Theorem 9(b), $n$-$\sigma X\backslash X$ is close to being extremally disconnected for if $V \in \tau(Y^{\#})$, $cl_{Y^{\#}}\backslash X$ is clopen in  $Y^{\#}\backslash X$.
\end{comment}

\subsection{Partitions of $n$-$\sigma X\setminus X$}

The theory of H-closed extensions starts with these two results that distinguish it from the theory of Hausdorff compactifications.  

Let $X$ be a Hausdorff space and $Z$ an H-closed extension of $X$.  Let $f:\kappa X \rightarrow Z$ be the continuous surjection such that $f(x) = x$ for $x \in X$. From \cite{PW} we have the following facts:

(a) Then $\mathcal{P} = \{f^{\leftarrow}(p):p \in Z\backslash X\}$ is a partition of $\sigma X \backslash
X$ into nonempty compact spaces.

(b)   If $\mathcal{P}$ is a partition of $\sigma X\backslash X$ into nonempty compact spaces, there is an H-closed extension $Z$ of $X$ such that if $f: \kappa X\rightarrow Z$ such that $f(x) = x$ for $x \in X$, then $\mathcal{P} = \{f^{\leftarrow}(p):p \in Z\backslash X\}$.

The question is whether this can be done with $n$-H-closed extensions. We give a partial answer to this problem considering the $n$-H-closed extensions that are Hausdorff except for $X$.  The following Lemma \ref{L4}, Theorems \ref{T10} and \ref{T11} show that it works in this very restricted setting.

\begin{lemma}\label{L4}
Let $n\in\omega$, $n\geq 2$, and $Z$ be an extension of  a space $X$.   

(a) If $U$ and $V$ are open in $Z$, then $U \subseteq o(U\cap X)$ and if $U \cap V = \varnothing$, then $o(U\cap X) \cap o(V\cap X) = \varnothing$.

(b) If $Z$ is Hausdorff except for $X$, then $Z^{\#}$ is Hausdorff except for $X$.

(c) If $Z$ is Hausdorff except for $X$ and $\mathcal{U}$ is an open ultrafilter on $X$, then $|a_Z\mathcal{U}\backslash a_X\mathcal{U}| \leq 1$.

(d) If $Z$ is Hausdorff except for $X$ and $Z$ is $n$-H-closed, then $n$-$\kappa X$ is Hausdorff except for $X$.
\end{lemma}
\proof The proof of (a) is straightforward.  For (b), let $p,q \in Z\backslash X$.  There are open sets $U,V$ such that $p \in U$, $q \in V$, and $U \cap V = \varnothing$.  Then $U \subseteq o(U\cap X)$,  $V \subseteq o(V\cap X)$, $p \in o(U\cap X)$, $q \in o(V\cap X)$, and $o(U\cap X) \cap o(V\cap X) = \varnothing$.  Similar proof for $p \in Z\backslash X$ and $x \in X$. (c). If $p,q \in a_Z\mathcal{U}\backslash a_X\mathcal{U}$, there are disjoint open sets $U,V$ in $Z$ such that $p \in U$ and $q \in V$.  Then $U\cap X, V\cap X \in \mathcal{U}$, a contradiction as $U \cap V = \varnothing$. (d). Let $p,q \in n$-$\kappa X \backslash X$.  Suppose $p, q \in a\mathcal{U}$, an open ultrafiler  $\mathcal{U}$ on $X$. Then $a_X\mathcal{U} \leq n-3$.  However, by (c), $a_X\mathcal{U} \geq n-2$, a contradiction.  So, $p \in a_X\mathcal{U}$ and $q \in a_X\mathcal{V}$ where $\mathcal{U}$ and $\mathcal{V}$ are distinct.  There are  $U \in \mathcal{U}$ and $V \in \mathcal{V}$ such that $U \cap V = \varnothing$.  Thus, points of $n$-$\kappa X \backslash X$ can be separated by disjoint open sets.  For a point $p \in n$-$\kappa X \backslash X$ and a point $x \in X$, $f(p)$ and $x$ can be separated by disjoint open sets. As $f$ is continuous, it follows that $p$ and $x$ can be separated by disjoint open sets.
\endproof

\begin{definition}\rm
A partition $\mathcal{P}$ of a subset of a space $X$ is said to be {\bf Hausdorff} if $A, B \in \mathcal{P}$, $A \ne B$, there are open sets $U,V$ in $X$ such that $A \subseteq U$, $B \subseteq V$, and $U \cap V = \varnothing$.
\end{definition}

\begin{theorem}\label{T10}
Let $n\in\omega$, $n\geq 2$, and $Z$ be an $n$-H-closed extension of a space $X$ that is Hausdorff except for $X$ and 
$f: n$-$\kappa X \rightarrow Z$ a continuous surjection such that $f(x) = x $. 
Then  $\mathcal{P} = \{f^{\leftarrow}(p): p \in Z\backslash X\}$ is a Hausdorff partition of compact subsets of  $n$-$\sigma X \backslash X$ with the property that each $f^{\leftarrow}(p)$ in $\mathcal{P}$ is  Hausdorff separated from  each point in $X$.
\end{theorem}
\proof
Fix $p \in Z\backslash X$, and let $K = f^{\leftarrow}(p)$.   First,  we show that $K$ is closed in $n$-$\sigma X$. Let $r \in n$-$\sigma X \backslash K$.  
As $Z$ is Hausdorff except for $X$, there are open sets $U, V \in \tau(Z)$ such that $f(r) \in V$ and $p \in U$ and $U\cap V = \varnothing$.  
As $f$ is $\theta$-continuous, there is an open set $T \in \tau(n$-$\sigma X$) such that $f[cl(T)] \subseteq cl(V)$ and $r \in T$.  As $ p \in U, f(r) \in f[cl(T)] \subseteq  cl(V)$ and $U \cap cl(V) = \varnothing$, $cl(T) \cap K = \varnothing$. A similar proof works when $r \in X$. This shows that $K$ is closed in $n$-$\sigma X$ and that $K$ is Hausdorff separated from each point in $X$. By Theorem \ref{9}(e), $K$ is compact as well as closed in $n$-$\sigma X$.
\endproof

\begin{theorem}\label{T11}
Let $n\in\omega$, $n\geq 2$, and $X$ be an $n$-Hausdorff space.  Let $\mathcal{P}$ be a Hausdorff partition of compact subsets of $n$-$\sigma X \backslash X$ that are closed in $n$-$\sigma X$ and Hausdorff separated from each point of $X$. There is an $n$-H-closed extension $Z$ of $X$ that is Hausdorff except for $X$ and a continuous surjection $f: n\text{-}\kappa X \rightarrow Z$ such that $f(x) = x $  and $\mathcal{P} = \{f^{\leftarrow}(p): p \in Z\backslash X\}$.
\end{theorem}
\proof
Let $\mathcal{P} = \{P_a:a \in A\}$ and $Z = X \cup A$ where $X \cap A = \varnothing$.  Define $f: n$-$\kappa X \rightarrow Z$ by $f(x) = x$ for $x \in X$ and if $p \in P_a$, $f(p) = a$.  The function $f$ is onto. 

We define $U \subseteq Z$ to be open if $U \cap X$ is open in $X$ and if $f(p) = a \in U$, $\{p\}\cup (U\cap X)$ is open in  $n$-$\kappa X$ for each $p \in P_a$.  Note that $Z$ is a simple extension of $X$, $f$ is continuous, and  if   $V$ is open in $X$ and  $\{p\}\cup (V)$ is open in  $n$-$\kappa X$ for each $p \in P_a$, $\{a\} \cup V$ is open in $Z$.  Conversely, if $\{a\} \cup V$ is open in $Z$, then  $\{p\}\cup (V)$ is open in  $n$-$\kappa X$ for each $p \in f^{\leftarrow}(a)$.

Next, we show that $Z$ is Hausdorff except for $X$:  Let $a,b \in Z\backslash X$.
There are open sets $U, V$ in $n$-$\sigma X$ such that $P_a \subseteq U$, $P_b \subseteq V$, and $U \cap V = \varnothing$.  Then $\{a\} \cup U\cap X$ and $\{b\} \cup V\cap X$ are open and disjoint in $Z$. A  similar proof shows that  $a \in A$ and $x \in X$ can be separated by disjoint disjoint open sets in $Z$ using the hypothesis that $P_a$ and $x$ can be separated by disjoint open sets in $n$-$\sigma X$. 

Next, we will show that $Z$ is $n$-Hausdorff. Let $A \subseteq Z$ such that $|A| = n$.  If $A \backslash X \ne \varnothing$, then there is a $a \in A \backslash X$.  Now $a$ can be separated by disjoint open sets from any point in $A \backslash X \cup\{a\}$ and any point in $X$.  The only remaining case is when $A \subseteq X$.  As $X$ is $n$-Hausdorff, there are open sets $\{U_a:a \in A\}$ in $X$ such that $a \in U_a$ and $\bigcap_{a\in A}U_a = \varnothing$.

Finally, we will show that $Z$ is $n$-H-closed.  Let $\mathcal{U}$ be an open ultrafilter on $Z$, by Theorem \ref{2}, it suffices to show that $|a_Z\mathcal{U}| = n-1$.  Let $\mathcal{V} = \{U\cap X: U \in \mathcal{U}\}$ is an open ultrafilter on $X$. Then   $|a_Z\mathcal{U}| = |a_Z\mathcal{V}| $.  
The goal is to show that $|a_{n\text{-}\kappa X}\mathcal{V}| = |a_Z\mathcal{V}|$, more precisely, that $|a_{n\text{-}\kappa X}\mathcal{V}\backslash a_X\mathcal{V} | = |a_Z\mathcal{V}\backslash a_X\mathcal{V}|$.  
Let   $ p \in a_{n\text{-}\kappa X}\mathcal{V}\backslash a_X\mathcal{V} $ and $f(p) = a \in Z\backslash X$.  As $ p \in a_{n\text{-}\kappa X}\mathcal{V}\backslash a_X\mathcal{V} $, $O^p_{{n\text{-}\kappa X}} = \mathcal{V}$. As $f(p) = a$ and $f$ is continuous, $O^a_Z \subseteq O^p_{{n\text{-}\kappa X}} = \mathcal{V}$.  Thus, $f(p) = a \in a_Z\mathcal{V}\backslash a_X\mathcal{V}$.
Conversely, suppose $a \in a_Z\mathcal{V}\backslash a_X\mathcal{V}$.  As $Z$ is $n$-Hausdorff, $|a_Z\mathcal{V}| \leq n-1$.  Thus  $|a_X\mathcal{V}| \leq n-2$. But, $|a_{n\text{-}\kappa X}\mathcal{V}| = n-1$. So, there is some $p \in a_{n\text{-}\kappa X}\mathcal{V}\backslash a_X\mathcal{V} $.
\endproof

\section*{Acknowledgment}
The third author was supported by the ``National Group for Algebraic and Geometric Structures, and their Applications'' (GNSAGA-INdAM).

\end{document}